\documentclass[11pt]{article}
\usepackage{amssymb,amsfonts,amsmath, psfrag,eepic}

\usepackage{pgf,tikz}
\usetikzlibrary{arrows.meta}

\usepackage{authblk}

\newtheorem{theo}{Theorem}

\newtheorem{prop}[theo]{Proposition}
\newtheorem{remark}{Remark}

\makeatletter \@addtoreset{equation}{section}
\@addtoreset{theo}{section} \makeatother

\def\qed{\hfill \rule{4pt}{7pt}}
\def\pf{\noindent {\it Proof.} }

\makeatletter
\newfont{\footsc}{cmcsc10 at 8truept}
\newfont{\footbf}{cmbx10 at 8truept}
\newfont{\footrm}{cmr10 at 10truept}

\makeatother 

\title{Involutions on Tip-Augmented Plane Trees for Leaf Interchanging}
\date{}

\author[1]{Laura L.M. Yang \footnote{\texttt{Email address:} laura.yang@ucf.edu}}
\author[2]{Dax T.X. Zhang \footnote{\texttt{Email address:} zhangtianxing6@tju.edu.cn}}
\affil[1]{\small Department of Mathematics, University of Central Florida, Orlando FL 32816, USA}
\affil[2]{\small College of Mathematical Science \& Institute of Mathematics and Interdisciplinary Sciences,\quad \quad  
Tianjin Normal University, Tianjin 300387, P.R. China}

\begin{document}
\maketitle

\begin{abstract} This paper constructs two involutions on tip-augmented plane trees, as defined by Donaghey, that interchange two distinct types of leaves while preserving all other leaves. These two involutions provide bijective explanations addressing a question posed by Dong, Du, Ji, and Zhang in their work.
\vskip 8pt

\end{abstract}

\section{Introduction}
Plane trees, also known as rooted ordered trees, are fundamental objects frequently studied in combinatorics, with many enumerative results appearing throughout the literature \cite{S}. For instance, the well-known Motkzin numbers \cite{D-S} enumerate plane trees in which the leftmost child of every interior vertex is a leaf.  Donaghey \cite{Don} referred to these structures as tip-augmented plane trees. 

Chen, Deutsch and Elizalde \cite{C-D-E} categorized the leaves of a plane tree into two types - old and young leaves - to refine the Narayana polynomials. More recently, Dong, Du, Ji, and Zhang \cite{D-D-J-Z} further refined this classification by dividing leaves into five distinct categories: singleton, elder twin, elder non-twin, second, and younger leaves. They derived generating functions for both refined Narayana polynomials and Motzkin polynomials, leading to the following symmetry between singleton and elder non-twin leaves in tip-augmented plane trees: 

\begin{theo} \label{main}
(Dong-Du-Ji-Zhang \cite{D-D-J-Z}) For $n\geq 2$, let $M_n(i, j, k, r, s)$ denote the number of tip-augmented plane trees with $n$ edges, $i$ singleton leaves, $j$ elder twin leaves, $k$ elder non-twin leaves, $r$ younger leaves and $s$ second leaves. Then 
$$M_n(i, j, k, r, s)=M_n(k, j, i, r, s).$$
\end{theo}
The objective of this paper is to construct two involutions on tip-augmented plane trees, providing bijective proofs for Theorem \ref{main}, as addressed in \cite{D-D-J-Z}. The first involution is constructed inductively on unlabelled tip-augmented plane trees. The second involution is constructed on labelled tip-augmented plane trees using the bijective algorithm provided by Chen \cite{Chen}. Both involutions interchange singleton leaves with elder non-twin leaves while leaving other leaves unchanged. 

To conclude this section, we review the definitions of the five categories of leaves in plane trees:
\begin{itemize}
    \item {\it Singleton leaf:} A leaf without any siblings.\\[-0.25in]
    \item {\it Elder twin leaf:} A leaf that is the leftmost child of its parent, where the second child of the same parent is also a leaf.\\[-0.25in]
    \item {\it Elder non-twin leaf:} A leaf that is the leftmost child of its parent, where the second child of the same parent is not a leaf.\\[-0.25in]
    \item {\it Second leaf:} A leaf that is the second child of its parent.\\[-0.25in]
    \item {\it Younger leaf:} A leaf that is neither the first nor the second child of its parent.
\end{itemize}

Figure~\ref{4 edges} illustrates four tip-augmented plane trees with four edges. The pair of numbers above each tree represents the numbers of singleton leaves (drawn with empty circles) and elder non-twin leaves (drawn with black filled circles), respectively. $\Phi$ is an involution constructed in Section 3.  
\vspace{-0.2in}
\begin{center}
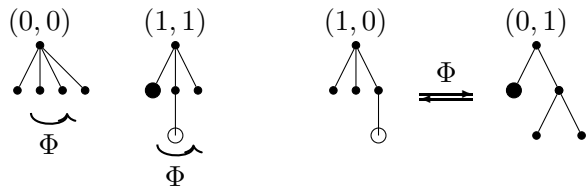
\begin{figure}[h]
\centering
\begin{tikzpicture}[scale=0.6]
\node at (0,1.5) {$(0,0)$};
\node[circle, draw, fill=black, inner sep=1pt] (r) at (0,1) {};
\foreach \i in {-1,...,2} {
\node[circle, draw, fill=black, inner sep=1pt](\i) at (\i/2,0) {};
\draw (r) -- (\i);
}

\draw[->, thick] (-0.2,-0.5) to[bend right=90] (0.6,-0.5);
\node at (0.2,-1.2) {$\Phi$}; 

\node at (3,1.5) {$(1,1)$};
\node[circle, draw, fill=black, inner sep=1pt] (r) at (3,1) {};
\foreach \i in {1} {
\node[circle, draw, fill=black, inner sep=2pt](\i) at (2+\i/2,0) {};
\draw (r) -- (\i);
}
\foreach \i in {2,3} {
\node[circle, draw, fill=black, inner sep=1pt](\i) at (2+\i/2,0) {};
\draw (r) -- (\i);
}
\draw (3,0)--(3,-1);
\node[circle, draw, inner sep=2pt] at (3,-1) {};

\draw[->, thick] (2.6,-1.2) to[bend right=90] (3.4,-1.2);
\node at (3,-1.9) {$\Phi$}; 

\node at (7,1.5) {$(1,0)$};
\node[circle, draw, fill=black, inner sep=1pt] (r) at (7,1) {};
\foreach \i in {1,...,3} {
\node[circle, draw, fill=black, inner sep=1pt](\i) at (6+\i/2,0) {};
\draw (r) -- (\i);
}
\draw (7.5,0)--(7.5,-1);
\node[circle, draw, inner sep=2pt] at (7.5,-1) {};

\node at (11,1.5) {$(0,1)$};
\node[circle, draw, fill=black, inner sep=1pt] (r) at (11,1) {};
\foreach \i in {1} {
\node[circle, draw, fill=black, inner sep=2pt](\i) at (10+\i/2,0) {};
\draw (r) -- (\i);
}
\foreach \i in {3} {
\node[circle, draw, fill=black, inner sep=1pt](\i) at (10+\i/2,0) {};
\draw (r) -- (\i);
}
\node[circle, draw, fill=black, inner sep=1pt] (s) at (11.5,0) {};
\foreach \i in {1,3} {
\node[circle, draw, fill=black, inner sep=1pt](\i) at (10.5+\i/2,-1) {};
\draw (s) -- (\i);
}
\node at (9,0.2){\vector(1,0){20}};
\node at (9,-0.2){\vector(-1,0){20}};
\node at (9,0.4){$\Phi$};
\end{tikzpicture}
    \caption{Involution $\Phi$ on tip-augmented plane trees with $4$ edges}
    \label{4 edges}
\end{figure}
\end{center}
\vspace{-0.3in}
Note that the categories of singleton, elder twin, and elder non-twin leaves refine the old leaves, while the categories of second and younger leaves refine the young leaves, as defined in \cite{C-D-E}.

\section{Classification for Tip-augmented Plane Trees}
Tip-augmented plane trees with $n\geq 2$ edges can be classified into one of the following four classes based on two criteria: (a) the number of edges of the subtree rooted at the second child of the root; (b) the existence of at least one child of the root that is a parent of a singleton leaf.  
\begin{itemize}
        \item [A1.] The second child of the root is a leaf, and no child of the root is a parent of a singleton leaf.
        \item [A2.] The second child of the root is a parent of a singleton leaf, and no other child of the root is a parent of a singleton leaf. 
        \item [B1.] The second child of the root is the parent of a subtree with at least two edges, and no child of the root is a parent of a singleton leaf.
        \item [B2.] Regardless of the second child of the root, at least one child of the root is a parent of a singleton leaf.  
    \end{itemize}

Tip-augmented plane trees are called {\it Type A} if they belong to class A1 or A2, and {\it Type B} if they belong to class B1 or B2. Figures~\ref{Type1} and \ref{Type2} illustrate the structures for Type A and Type B, respectively. 

It is worth noting that each subtree of the root, represented by a triangle, is either a trivial tree with one vertex or a tip-augmented plane tree with at least two edges. A tree depicted with an ellipse contains at least two edges. In particular, in class B2, vertex $u$ is the rightmost child of the root that has only one child which is a singleton leaf. Hence, $B$ is a tree rooted at~$r$,  positioned to the left of the edge $ru$.  

\begin{center}
\begin{figure}[ht]
\centering
\begin{tikzpicture}[scale=0.8]
\node at (4,1.5) {\small{Class A1}};
\node[circle, draw, fill=black, inner sep=1pt] (r) at (4,1) {};
\node[circle, draw, fill=black, inner sep=1pt] (s) at (4,0) {};
\node[circle, draw, fill=black, inner sep=1pt] (t) at (3,0) {};
\draw (r) --(s);
\draw(r)--(t);

\foreach \i in {5,6} {
\node[circle, draw, fill=black, inner sep=1pt](\i) at (\i,0) {};
\draw (r) -- (\i)--(\i-0.3, -1)--(\i+0.3,-1)--(\i);
}
\draw (r) -- (7.2, 0) --(6.9,-1) -- (7.5, -1) -- (7.2, 0);

\node[circle, draw, fill=black, inner sep=1pt] at (7.2,0) {};
\node at (5, -1.3) {\scriptsize{$T_1$}};
\node at (6, -1.3) {\scriptsize{$T_2$}};
\node at (6.5, 0) {$\cdots$};
\node at (7.3, -1.3) {\scriptsize{$T_k$}};

\node at (12, 1.5) {\small{Class A2}};
\node[circle, draw, fill=black, inner sep=1pt] (r) at (12,1) {};
\node[circle, draw, fill=black, inner sep=1pt] (s) at (12,0) {};
\node[circle, draw, fill=black, inner sep=2pt] (t) at (11,0) {};
\node[circle, draw, inner sep=2pt] at (12,-1) {};
\draw(r)--(s);
\draw(r)--(t);
\draw(r)--(12,-1);

\node at (12.3,0) {$u$};
\node at (12, -1.4) {$v$};
\node at (11, -0.4) {$w$};

\foreach \i in {13,14} {
\node[circle, draw, fill=black, inner sep=1pt](\i) at (\i,0) {};
\draw (r) -- (\i)--(\i-0.3, -1)--(\i+0.3,-1)--(\i);
}
\draw (r) -- (15.2, 0) --(14.9,-1) -- (15.6, -1) -- (15.2, 0);
\node[circle, draw, fill=black, inner sep=1pt] at (15.2,0) {};
\node at (13, -1.3) {\scriptsize{$T_1$}};
\node at (14, -1.3) {\scriptsize{$T_2$}};
\node at (14.5, 0) {$\cdots$};
\node at (15.3, -1.3) {\scriptsize{$T_k$}};
\end{tikzpicture}
    \caption{Tip-augmented plane trees of Type A}
    \label{Type1}
\end{figure}
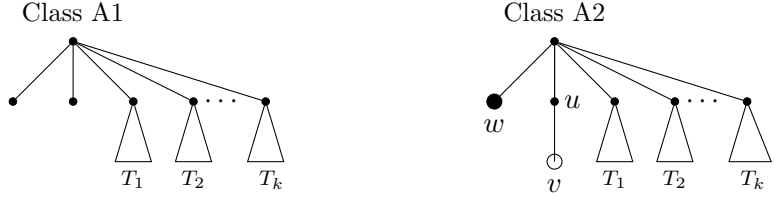
\end{center}
\vspace{-0.4in}
\begin{center}
\begin{figure}[ht]
\centering
\begin{tikzpicture}[scale=0.8]
\node at (4,1.8) {\small{Class B1}};
\node[circle, draw, fill=black, inner sep=1pt, label=$r$] (r) at (4,1) {};
\node[circle, draw, fill=black, inner sep=1pt] (u) at (4,0) {};
\node[circle, draw, fill=black, inner sep=2pt] (v) at (3,0) {};
\draw(r)--(u);
\draw(r)--(v);
\node at (3,-0.3) {$w$};
\node at (4.4,0.1) {$u$};
\draw[dashed] (4,0) ellipse (0.2cm and 0.2cm);

\draw (4,-0.5) ellipse (0.3cm and 0.5cm);
\node at (4, -1.3) {\scriptsize{$A$}};

\foreach \i in {5,6} {
\node[circle, draw, fill=black, inner sep=1pt](\i) at (\i,0) {};
\draw (r) -- (\i)--(\i-0.3, -1)--(\i+0.3,-1)--(\i);
}
\draw (r) -- (7.2, 0) --(6.9,-1) -- (7.5, -1) -- (7.2, 0);
\node[circle, draw, fill=black, inner sep=1pt] at (7.2,0) {};
\node at (5, -1.3) {\scriptsize{$T_1$}};
\node at (6, -1.3) {\scriptsize{$T_2$}};
\node at (6.4, 0) {$\cdots$};
\node at (7.2, -1.3) {\scriptsize{$T_k$}};
\node at (12,1.8) {\small{Class B2}};
\node[circle, draw, fill=black, inner sep=1pt, label=$r$] (r) at (12,1) {};
\node[circle, draw, fill=black, inner sep=1pt] (u) at (12,0) {};
\node[circle, draw, inner sep=2pt] (s) at (12,-1) {};
\draw(r)--(s);
\node at (12.4,0) {$u$};
\node at (12,-1.3) {$v$};
\draw[dashed] (12,0) ellipse (0.2cm and 0.2cm);

\foreach \i in {13,14} {
\node[circle, draw, fill=black, inner sep=1pt](\i) at (\i,0) {};
\draw (r) -- (\i)--(\i-0.3, -1)--(\i+0.3,-1)--(\i);
}
\draw (r) -- (15.2, 0) --(14.9,-1) -- (15.5, -1) -- (15.2, 0);

\node[circle, draw, fill=black, inner sep=1pt] at (15.2,0) {};
\node at (13, -1.3) {\scriptsize{$T'_1$}};
\node at (14.5, 0) {$\cdots$};
\node at (14, -1.3) {\scriptsize{$T'_2$}};
\node at (15.3, -1.3) {\scriptsize{$T'_k$}};

\draw[rotate=45] 
(8.25,-7.75) ellipse (0.9cm and 0.3cm);
\node at (11, -0.5) {\scriptsize{$B$}};
\end{tikzpicture}
    \caption{Tip-augmented plane trees of Type B}
    \label{Type2}
\end{figure}
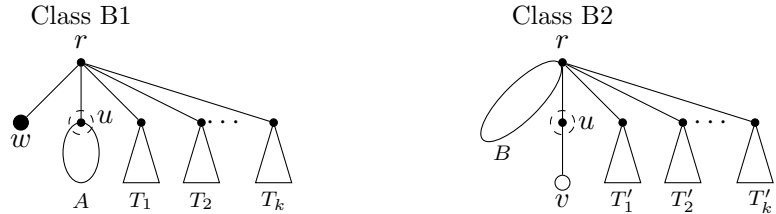
\end{center}

\vspace{-0.3in}  
\begin{prop}\label{Prop} Assume that $T$ is a tip-augmented plane tree with $n\geq 2$ edges. There are four cases for counting singleton and elder non-twin leaves, depending on the specific class of $T$.  \\

\noindent {Class A1:} The number of singleton (elder non-twin) leaves equals the sum of the numbers of singleton (elder non-twin) leaves in all subtrees $T_i$.\\

\noindent{Class A2:} The number of singleton (elder non-twin) leaves equals one plus the sum of the numbers of singleton (elder non-twin) leaves in the subtrees $T_i$.\\ 

\noindent {Class B1:} The number of singleton leaves equals the sum of the numbers of singleton leaves in tree $A$ and all subtrees $T_i$. The number of elder non-twin leaves equals one plus the sum of the numbers of elder non-twin leaves in tree $A$ and all subtrees $T_i$.\\

\noindent{Class B2:} The number of singleton leaves equals one plus the sum of the numbers of singleton leaves in tree $B$ and all subtrees $T'_i$. The number of elder non-twin leaves equals the sum of the numbers of elder non-twin leaves in tree $A'$ and all subtrees $T'_i$.
\end{prop}
\pf For any tip-augmented plane tree $T$ with $n\geq 2$ edges, every singleton (elder non-twin) leaf in each subtree of the root is also a singleton (elder non-twin) leaf in tree~$T$.~\qed

\section{An Involution $\Phi$ on Tip-augmented Plane Trees}
Based on the property that any subtree rooted at an interior vertex of a tip-augmented plane tree is also a tip-augmented plane tree, an involution $\Phi$ can be constructed inductively on the number of edges. 
\begin{theo}
There is an involution $\Phi$ on the set of tip-augmented plane trees with $n$ edges that interchanges elder non-twin leaves and singleton leaves while preserving elder twin leaves, younger leaves, and second leaves.
\end{theo}
\pf (Induction on $n$) Define $\Phi(T)=T$ for $n=0, 1, 2$. It is straightforward to verify that, in these cases, both the number of singleton leaves and the number of elder non-twin leaves are equal to 0. Now, let $T$ be a tip-augmented plane tree with $n$ edges for $n\geq 3$. The involution $\phi(T)$ is constructed recursively based on the class to which $T$ belongs. 

Assume that $T$ is a tree of Type A with a root of degree $k+2$, as illustrated in Figure~\ref{Type1}. The second subtree of the root is either a trivial tree or a tree with one edge, and each $T_i$ is either a trivial tree or a tip-augmented plane tree with at least two edges for $i=1, \ldots, k$. The image $\Phi(T)$ is constructed by replacing each $T_i$ with $\Phi(T_i)$, while the first two subtrees of the root remain unchanged. Thus, $\Phi(T)$ belongs to the same class as $T$, and the type of each leaf in the subtrees $\Phi(T_i)$ remains unchanged in $\Phi(T)$. By induction, the singleton leaves of each $T_i$ are converted into elder non-twin leaves in $\Phi(T_i)$, while the type of all other leaves remains unchanged. Consequently, the singleton leaves and elder non-twin leaves of $T$ are swapped in $\Phi(T)$ while all other leaves are preserved, assuming that the vertices $w$ and $v$ are swapped in class A2.   

Assume that $T$ is a tree of type B with a root $r$ of degree $k+2$. If $T$ is a tree of class B1, whose structure is illustrated in Figure~\ref{Type2}, then the first child $w$ is an elder non-twin leaf, the second child $u$ is the root of tree $A$ which is a tip-augmented plane tree with at least two edges, and each $T_i$ is either a trivial tree or a tip-augmented plane tree with at least two edges. The image $\Phi(T)$ is constructed as follows: place $\Phi(A)$ to the left of the edge $ru$ by identifying $r$ as the root, add a vertex $v$ as a child of $u$, and set $\Phi(T_i)$ as the subtrees of the root lying to the right of edge $ru$ in the same order as in $T$. 

Comparing $\Phi(T)$ with $T$, the involution $\Phi$ converts the elder non-twin leaf $w$ in $T$ into a singleton leaf $v$ in $\Phi(T)$, and the second child $u$ in $T$ becomes the rightmost child of the root in $\Phi(T)$, where $u$ is the parent of the singleton leaf $v$. Consequently, $\Phi(T)$ is a tree of class B2.  Moreover, the type of each leaf in $\Phi(A)$ and $\Phi(T_i)$ remains unchanged in $\Phi(T)$. By induction, the singleton leaves and elder non-twin leaves of $T$ are swapped in $\Phi(T)$, while all other leaves retain their types. If $T$ is a tree of class B2, the above construction can be reversed to obtain $\Phi(T)$, completing the involution.
\qed

\begin{remark}
The involution $\Phi$ maps trees of type B between class B1 and class B2, while trees of type A remain in their original class.
\end{remark}

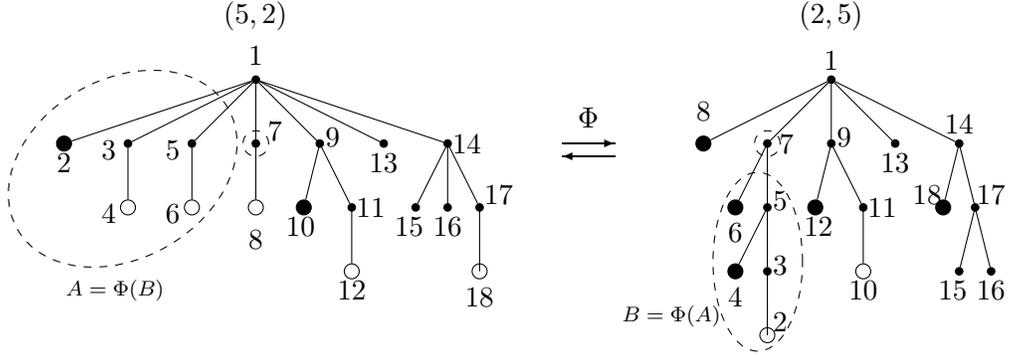
\begin{figure}[ht]
\centering
\begin{tikzpicture}[scale=0.85]
\node at (0,2) {$(5,2)$};
\node[circle, draw, fill=black, inner sep=1pt, label=1] (r) at (0,1) {};
\foreach \i in {-3} {
\node[circle, draw, fill=black, inner sep=2pt](\i) at (\i,0) {};
\draw (r) -- (\i);
}
\foreach \i in {-2,...,3} {
\node[circle, draw, fill=black, inner sep=1pt](\i) at (\i,0) {};
\draw (r) -- (\i);
}
\foreach \i in {0} {
\node[circle, draw, inner sep=2pt](\i) at (-2+\i,-1) {};
\draw (-2,0) -- (\i);
}
\foreach \i in {0} {
\node[circle, draw, inner sep=2pt](\i) at (-1+\i,-1) {};
\draw (-1,0) -- (\i);
}
\draw[rotate=30,dashed] (-2,0.7) ellipse (1.9cm and 1.4cm);
\node at (-2.2, -2.3){\scriptsize{$A=\Phi(B)$}};
\node at (-3,-0.3) {$2$};
\node at (-2.3,-0.1) {$3$};
\node at (-2.3,-1.1) {$4$};
\node at (-1.3,-0.1) {$5$};
\node at (-1.3,-1.1) {$6$};

\foreach \i in {0} {
\node[circle, draw, inner sep=2pt](\i) at (\i,-1) {};
\draw (0,0) -- (\i);
}
\node at (0.3,0.2) {$7$};
\node at (0,-1.5) {$8$};
\draw[dashed] (0,0) ellipse (0.2cm and 0.2cm);

\foreach \i in {-1} {
\node[circle, draw, fill=black, inner sep=2pt](\i) at (1+\i/4,-1) {};
\draw (1,0) -- (\i);
}
\foreach \i in {1} {
\node[circle, draw, fill=black, inner sep=1pt](\i) at (1+\i/2,-1) {};
\draw (1,0) -- (\i);
}

\foreach \i in {0} {
\node[circle, draw, inner sep=2pt](\i) at (1.5+\i,-2) {};
\draw (1.5,-1) -- (\i);
}
\node at (1.2,0.1) {$9$};
\node at (0.7,-1.3) {$10$};
\node at (1.8,-1) {$11$};
\node at (1.5,-2.3) {$12$};
\node at (2,-0.3) {$13$};

\foreach \i in {-1,...,1} {
\node[circle, draw, fill=black, inner sep=1pt](\i) at (3+\i/2,-1) {};
\draw (3,0) -- (\i);
}
\node[circle, draw,  inner sep=2pt] at (3.5,-2) {};
\draw (3.5,-1) -- (3.5,-2);
\node at (3.3,0) {$14$};
\node at (2.4,-1.3) {$15$};
\node at (3,-1.3) {$16$};
\node at (3.8,-0.8) {$17$};
\node at (3.5,-2.4) {$18$};

\node at (5.2,0.2){\vector(1,0){20}};
\node at (5.2,-0.2){\vector(-1,0){20}};
\node at (5.2,0.4){$\Phi$};
\node at (9,2) {$(2,5)$};
\node[circle, draw, fill=black, inner sep=1pt] (r) at (9,1) {};
\foreach \i in {-2} {
\node[circle, draw, fill=black, inner sep=2pt](\i) at (9+\i,0) {};
\draw (r) -- (\i);
}
\foreach \i in {-1,...,2} {
\node[circle, draw, fill=black, inner sep=1pt](\i) at (9+\i,0) {};
\draw (r) -- (\i);
}
\node at (9,1.3) {$1$};
\node at (7,0.5) {$8$};

\foreach \i in {-1} {
\node[circle, draw, fill=black, inner sep=2pt](\i) at (8+\i/2,-1) {};
\draw (8,0) -- (\i);
}
\foreach \i in {0} {
\node[circle, draw, fill=black, inner sep=1pt](\i) at (8+\i/2,-1) {};
\draw (8,0) -- (\i);
}
\foreach \i in {-1} {
\node[circle, draw, fill=black, inner sep=2pt](\i) at (8.+\i/2,-2) {};
\draw (8,-1) -- (\i);
}
\foreach \i in {0} {
\node[circle, draw, fill=black, inner sep=1pt](\i) at (8.+\i/2,-2) {};
\draw (8,-1) -- (\i);
}
\node[circle, draw, inner sep=2pt]at (8,-3) {};
\draw (8,-2) -- (8,-3);

\draw[dashed] (8,0) circle (0.2cm and 0.2cm);
\draw[dashed] (7.85,-1.85) ellipse (0.7cm and 1.4cm);
\node at (6.5, -2.7){\scriptsize{$B=\Phi(A)$}};

\node at (8.3,0) {$7$};
\node at (7.5,-1.4) {$6$};
\node at (7.5,-2.4) {$4$};
\node at (8.2,-0.9) {$5$};
\node at (8.2,-1.9) {$3$};
\node at (8.2,-2.8) {$2$};
\foreach \i in {-1} {
\node[circle, draw, fill=black, inner sep=2pt](\i) at (9+\i/4,-1) {};
\draw (9,0) -- (\i);
}
\foreach \i in {1} {
\node[circle, draw, fill=black, inner sep=1pt](\i) at (9+\i/2,-1) {};
\draw (9,0) -- (\i);
}

\foreach \i in {0} {
\node[circle, draw, inner sep=2pt](\i) at (9.5+\i,-2) {};
\draw (9.5,-1) -- (\i);
}
\node at (9.2,0.1) {$9$};
\node at (8.8,-1.3) {$12$};
\node at (9.8,-1) {$11$};
\node at (9.5,-2.3) {$10$};
\node at (10,-0.3) {$13$};

\foreach \i in {-1} {
\node[circle, draw, fill=black, inner sep=2pt](\i) at (11+\i/4,-1) {};
\draw (11,0) -- (\i);
}
\foreach \i in {1} {
\node[circle, draw, fill=black, inner sep=1pt](\i) at (11+\i/4,-1) {};
\draw (11,0) -- (\i);
}
\foreach \i in {-1,1} {
\node[circle, draw, fill=black, inner sep=1pt](\i) at (11.25+\i/4,-2) {};
\draw (11.25,-1) -- (\i);
}
\node at (11,0.3) {$14$};
\node at (10.5,-0.8) {$18$};
\node at (11.5,-0.8) {$17$};
\node at (10.9,-2.3) {$15$};
\node at (11.5,-2.3) {$16$};
\end{tikzpicture}
    \caption{An example of involution $\Phi$}
    \label{n=17}
\end{figure}

Figures~\ref{4 edges} and \ref{5 edges} illustrate the involution $\Phi$ applied to tip-augmented plane tree with four and five edges, respectively. Figure~\ref{n=17} demonstrates how each vertex in $T$ is mapped to the corresponding vertex in $\Phi(T)$, preserving the same label.
\begin{center}
\begin{figure}[ht]
\centering
\begin{tikzpicture}[scale=0.6]
\node at (0,1.5) {$(0, 0)$};
\node[circle, draw, fill=black, inner sep=1pt] (r) at (0,1) {};
\foreach \i in {-2,...,2} {
\node[circle, draw, fill=black, inner sep=1pt](\i) at (\i/2,0) {};
\draw (r) -- (\i);
}

\draw[->, thick] (-0.6,-0.5) to[bend right=90] (0.4,-0.5);
\node at (0,-1) {$\Phi$}; 

\node at (2.5,1.5) {$(1, 1)$};
\node[circle, draw, fill=black, inner sep=1pt] (r) at (2.5,1) {};

\foreach \i in {-1} {
\node[circle, draw, fill=black, inner sep=2pt](\i) at (2.5+\i/2,0) {};
\draw (r) -- (\i);
}
\foreach \i in {0,...,2} {
\node[circle, draw, fill=black, inner sep=1pt](\i) at (2.5+\i/2,0) {};
\draw (r) -- (\i);
}
\foreach \i in {0} {
\node[circle, draw, inner sep=2pt](\i) at (2.5+\i/2,-1) {};
\draw (2.5,0) -- (\i);
}
\draw[->, thick] (2.1,-1.2) to[bend right=90] (2.9,-1.2);
\node at (2.5,-1.7) {$\Phi$}; 

\node at (5,1.5) {$(0, 0)$};
\node[circle, draw, fill=black, inner sep=1pt] (r) at (5,1) {};
\foreach \i in {-1,...,1} {
\node[circle, draw, fill=black, inner sep=1pt](\i) at (5+\i/2,0) {};
\draw (r) -- (\i);
}
\foreach \i in {-1,1} {
\node[circle, draw, fill=black, inner sep=1pt](\i) at (5.5+\i/4,-1) {};
\draw (5.5,0) -- (\i);
}
\draw[->, thick] (4.6,-1.2) to[bend right=90] (5.4,-1.2);
\node at (5,-1.7) {$\Phi$}; 

\node at (7.5,1.5) {$(1, 0)$};
\node[circle, draw, fill=black, inner sep=1pt] (r) at (7.5,1) {};
\foreach \i in {-1,...,2} {
\node[circle, draw, fill=black, inner sep=1pt](\i) at (7.5+\i/2,0) {};
\draw (r) -- (\i);
}
\draw (8,0)--(8,-1);
\node[circle, draw, inner sep=2pt] at (8,-1) {};

\node at (9.5,0.2){\vector(1,0){15}};
\node at (9.5,-0.2){\vector(-1,0){15}};
\node at (9.5,0.4){$\Phi$};
\node at (11,1.5) {$(0, 1)$};
\node[circle, draw, fill=black, inner sep=1pt] (r) at (11,1) {};
\foreach \i in {-1} {
\node[circle, draw, fill=black, inner sep=2pt](\i) at (11+\i/2,0) {};
\draw (r) -- (\i);
}
\foreach \i in {0,1} {
\node[circle, draw, fill=black, inner sep=1pt](\i) at (11+\i/2,0) {};
\draw (r) -- (\i);
}
\foreach \i in {-1,1} {
\node[circle, draw, fill=black, inner sep=1pt](\i) at (11+\i/4,-1) {};
\draw (11,0) -- (\i);
}
\node at (13.5,1.5) {$(1, 0)$};
\node[circle, draw, fill=black, inner sep=1pt] (r) at (13.5,1) {};
\foreach \i in {-1,...,2} {
\node[circle, draw, fill=black, inner sep=1pt](\i) at (13.5+\i/2,0) {};
\draw (r) -- (\i);
}
\foreach \i in {0} {
\node[circle, draw, inner sep=2pt](\i) at (14.5,-1) {};
\draw (14.5,0) -- (\i);
}
\node at (15.5,0.2){\vector(1,0){15}};
\node at (15.5,-0.2){\vector(-1,0){15}};
\node at (15.5,0.4){$\Phi$};
\node at (17,1.5) {$(0, 1)$};
\node[circle, draw, fill=black, inner sep=1pt] (r) at (17,1) {};
\foreach \i in {-1} {
\node[circle, draw, fill=black, inner sep=2pt](\i) at (17+\i/2,0) {};
\draw (r) -- (\i);
}
\foreach \i in {0} {
\node[circle, draw, fill=black, inner sep=1pt](\i) at (17+\i/2,0) {};
\draw (r) -- (\i);
}
\node[circle, draw, fill=black, inner sep=1pt] (s) at (17,0) {};
\foreach \i in {-1,...,1} {
\node[circle, draw, fill=black, inner sep=1pt](\i) at (17+\i/2,-1) {};
\draw (s) -- (\i);
}
\node at (19.5,1.5) {$(2, 1)$};
\node[circle, draw, fill=black, inner sep=1pt] (r) at (19.5,1) {};
\foreach \i in {-1} {
\node[circle, draw, fill=black, inner sep=2pt](\i) at (19.5+\i/2,0) {};
\draw (r) -- (\i);
}
\foreach \i in {0,...,1} {
\node[circle, draw, fill=black, inner sep=1pt](\i) at (19.5+\i/2,0) {};
\draw (r) -- (\i);
}
\draw (19.5,0)--(19.5, -1);
\node[circle, draw, inner sep=2pt] (s) at (19.5,-1) {};
\foreach \i in {0} {
\node[circle, draw, inner sep=2pt](\i) at (20,-1) {};
\draw (20,0) -- (\i);
}
\node at (21,0.2){\vector(1,0){15}};
\node at (21,-0.2){\vector(-1,0){15}};
\node at (21,0.4){$\Phi$};
\node at (22.5,1.5) {$(1, 2)$};
\node[circle, draw, fill=black, inner sep=1pt] (r) at (22.5,1) {};
\foreach \i in {-1} {
\node[circle, draw, fill=black, inner sep=2pt](\i) at (22.5+\i/4,0) {};
\draw (r) -- (\i);
}
\foreach \i in {1} {
\node[circle, draw, fill=black, inner sep=1pt](\i) at (22.5+\i/4,0) {};
\draw (r) -- (\i);
}
\foreach \i in {-1} {
\node[circle, draw, fill=black, inner sep=2pt](\i) at (22.75+\i/4,-1) {};
\draw (22.75,0) -- (\i);
}
\foreach \i in {1} {
\node[circle, draw, fill=black, inner sep=1pt](\i) at (22.75+\i/4,-1) {};
\draw (22.75,0) -- (\i);
}
\foreach \i in {0} {
\node[circle, draw, inner sep=2pt](\i) at (23,-2) {};
\draw (23, -1) -- (\i);
}
\end{tikzpicture}
    \caption{Involution $\Phi$ on tip-augmented plane trees with $5$ edges}
    \label{5 edges}
\end{figure}
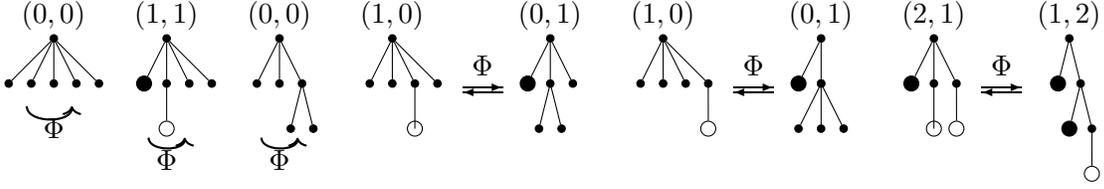
\end{center}
\vspace{-0.3in}

\section{An Involution $\Psi$ on Labelled Tip-augmented Plane Trees}
In this section, we present an involution on labelled tip-augmented plane trees that interchanges elder non-twin leaves and singleton leaves while preserving elder twin leaves, younger leaves, and second leaves, using Chen's bijective algorithm for labelled plane trees in \cite{Chen}.

Chen's bijective algorithm decomposes a labeled plane tree on $\{1, 2, \ldots, n+1\}$ into a set $F$ of $n$ matches with labels $\{1, \ldots, n, n+1, (n+2)^*, \dots, (2n)^*\}$, where a match is a rooted tree with two vertices. The reverse procedure of this decomposition algorithm is the following merging algorithm.  We start with a set $F$ of matches on $\{1, \ldots, n+1, (n+2)^*, \ldots, (2n)^*\}$. A vertex labelled by a mark $*$ is called a marked vertex.\\[-3mm]

\noindent (1) Find the tree $T$ with the smallest root in which no
vertex is marked. Let $i$ be the root of
$T$.\\[-3mm]

\noindent (2) Find the tree $T^*$ in $F$ that contains the smallest marked
vertex. Let $j^*$ be this marked vertex.\\[-3mm]

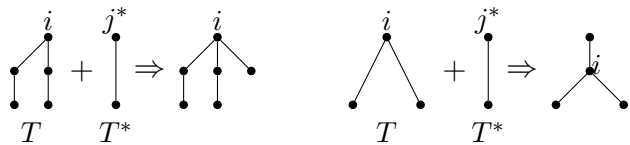
\begin{figure}[ht]
\centering
\begin{tikzpicture}[scale=0.45]
\node at (0,1.5) {$i$};
\node[circle, draw, fill=black, inner sep=1pt] (r1) at (0,1) {};

\foreach \i in {-1, 0} 
{
\node[circle, draw, fill=black, inner sep=1pt](\i) at (\i,0) {};
\node[circle, draw, fill=black, inner sep=1pt] at (\i,-1) {};
\draw (r1) -- (\i);
\draw (\i) -- (\i,-1);
}

\node at (-0.5, -1.8) {$T$};

\node at (1,0){\large\bf$+$};
\node at (2,1.5){$j^*$};
\node[circle, draw, fill=black, inner sep=1pt] (r2) at (2,1) {};
\node[circle, draw, fill=black, inner sep=1pt](s) at (2,-1) {};
\draw (r2) -- (s);
\node at (2, -1.8) {$T^*$};

\node at (3,0){\large\bf$\Rightarrow$};

\node at (5,1.5) {$i$};
\node[circle, draw, fill=black, inner sep=1pt] (r3) at (5,1) {};
\foreach \i in {4, 5} {
\node[circle, draw, fill=black, inner sep=1pt](\i) at (\i,0) {};
\node[circle, draw, fill=black, inner sep=1pt] at (\i,-1) {};
\draw (r3) -- (\i);
\draw (\i) -- (\i,-1);
}

\node[circle, draw, fill=black, inner sep=1pt] at (6,0) {};
\draw (r3) -- (6,0);

\node at (10,1.5) {$i$};
\node[circle, draw, fill=black, inner sep=1pt] (r1) at (10,1) {};
\foreach \i in {9, 11} {
\node[circle, draw, fill=black, inner sep=1pt](\i) at (\i,-1) {};
\draw (r1) -- (\i);
}
\node at (10, -1.8) {$T$};

\node at (12,0){\large\bf$+$};
\node at (13,1.5){$j^*$};
\node[circle, draw, fill=black, inner sep=1pt] (r2) at (13,1) {};
\node[circle, draw, fill=black, inner sep=1pt](t) at (13,-1) {};
\draw (r2) -- (t);
\node at (13, -1.8) {$T^*$};

\node at (14,0){\large\bf$\Rightarrow$};

\node at (16.2,0.2) {$i$};
\node[circle, draw, fill=black, inner sep=1pt] (r5) at (16,1) {};
\node[circle, draw, fill=black, inner sep=1pt] at (16,0) {};
\draw (r5) -- (16,0);

\foreach \i in {15, 17} {
\node[circle, draw, fill=black, inner sep=1pt](\i) at (\i,-1) {};
\draw (16,0) -- (\i);
}
\end{tikzpicture}
\caption{A horizontal merge and vertical merge}\label{merge}
\end{figure}

\noindent (3) If $j^*$ is the root of $T^*$, then merge $T$ and $T^*$ by
identifying $i$ and $j^*$, keep $i$ as the new vertex, and put the
subtrees of $T^{*}$ to the right of $T$. The operation is called a
\emph{horizontal merge}. If $j^*$ is a leaf of $T^*$, then replace
$j^*$ with $T$ in $T^*$. This operation is called a \emph{vertical
merge}. See Figure~\ref{merge}.\\[-3mm]

\noindent (4) Repeat the above procedure until $F$ becomes a
labelled tree.\\[-3mm]

For any set $F$ of $n$ matches labelled by $\{1, \ldots, n+1, (n+2)^*, \ldots, (2n)^*\}$, there are four types of matches when decomposing a plane tree $T$ into matches:
\begin{itemize}
    \item[i.] A match consists of two unmarked vertices,one of which is  a leaf corresponding to either an elder leaf or a singleton leaf in $T$.
    \item[ii.] A match consists of a marked root and an unmarked leaf corresponding to a young leaf in $T$.
    \item[iii.] A match consists of a marked leaf and an unmarked root corresponding to an old interior vertex in $T$. 
    \item[iv.] A match consists of two marked vertices, with the root corresponding to a young interior vertex in $T$.
\end{itemize}
Based on the observation above, the authors \cite{C-S-Y} provide an involution on labelled plane trees by turning a special match upside down. This match consists of a marked vertex and an unmarked vertex with the minimum value. The upside-down operation interchanges a match of type ii match with a  match of type iii that reverse the parity of young leaves. This gives a bijective explanation for Corollary 1.3 in \cite{D-D-J-Z}. 

When Chen's bijective algorithm is applied exclusively to the set of labelled tip-augmented plane trees, we can easily verify the following properties based on the definition of tip-augmented plane tree, where the leftmost child of each interior vertex is a leaf.

\begin{prop}
The following properties hold when applying Chen's bijective algorithm to tip-augmented planed trees.
\begin{itemize}
 \item[a)] No matches of type \textrm{iii}  are obtained when applying the decomposition algorithm to any tip-augmented plane tree.
 \item[b)] An elder non-twin leaf arises only when horizontally merging a match of type \textrm{i} with a match of type \textrm{iv}. In other words, the label of the root in a match of type \textrm{iv} is smaller than the label of its leaf. 
 \item[c)] A singleton leaf arises only when vertically merging a match of type {\textrm i} with a match of type \textrm{iv}. In other words, the label of  the root in a match of iv is greater than the label of its leaf.  
\end{itemize}
\end{prop}
\begin{theo}
    There exists an involution $\Psi$ on the set of labelled tip-augmented plane trees that interchanges elder non-twin leaves and singleton leaves while preserving elder twin leaves, younger leaves, and second leaves. 
\end{theo}
\pf The involution $\Psi$ is defined on the set of labelled tip-augmented plane trees whose match decomposition excludes any matches of type \textrm{iii}. Given a tip-augmented plane tree $T$, we first decompose it into matches, then turn all matches of type \textrm{iv} upside down. The image $\Psi(T)$ is obtained by applying the merging algorithm to the set of modified matches. 
\qed

\begin{figure}[ht]
\centering
\begin{tikzpicture}[scale=0.8]
\node at (0,2) {$(1,3)$};
\node[circle, draw, fill=black, inner sep=1pt] (r) at (0,1) {};
\node at (0,1.3) {$10$};

\node[circle, draw, fill=black, inner sep=2pt] at (-2,0) {};
\foreach \i in {-2,...,1} {
\node[circle, draw, fill=black, inner sep=1pt](\i) at (\i,0) {};
\draw (r) -- (\i);
}
\node at (-2.3,-0.1) {$3$};
\node at (-1.3,-0.1) {$1$};
\node at (0.3,-0.1) {$9$};
\node at (1.3,-0.1) {$7$};

\foreach \i in {-1, 1} {
\node[circle, draw, fill=black, inner sep=1pt](\i) at (-1+\i/2,-1) {};
\draw (-1,0) -- (\i);
}
\node at (-1.5,-1.3) {$8$};
\node at (-0.5,-1.3) {$5$};

\foreach \i in {-1} {
\node[circle, draw, fill=black, inner sep=2pt](\i) at (1+\i/4,-1) {};
\draw (1,0) -- (\i);
}
\foreach \i in {1} {
\node[circle, draw, fill=black, inner sep=1pt](\i) at (1+\i/2,-1) {};
\draw (1,0) -- (\i);
}

\foreach \i in {-1} {
\node[circle, draw, fill=black, inner sep=2pt](\i) at (1.5+\i/4,-2) {};
\draw (1.5,-1) -- (\i);
}

\foreach \i in {1} {
\node[circle, draw, fill=black, inner sep=1pt](\i) at (1.5+\i/2,-2) {};
\draw (1.5,-1) -- (\i);
}
\node[circle, draw, inner sep=2pt] at (2,-3) {};
\draw (2,-2) -- (2, -3);
\node at (0.7,-1.3) {$6$};
\node at (1.8,-1) {$2$};
\node at (1,-2) {$4$};
\node at (2.4,-2) {$12$};
\node at (2.4,-3) {$11$};
\node at (-0.1,-4){\vector(0,1){20}};
\node at (0.1,-4){\vector(0,-1){20}};
\foreach \i in {-5,...,5} {
\node[circle, draw, fill=black, inner sep=1pt](\i) at (\i/1.5,-5.5) {};
\node[circle, draw, fill=black, inner sep=1pt](-\i) at (\i/1.5,-6.5) {};
\draw (\i) -- (-\i);
}
\node at (-3.3,-5.1){$1$}; \node at (-3.3,-6.9){$8$};
\node at (-2.7,-5.1){$13^*$}; \node at (-2.7,-6.9){$5$};
\node at (-2,-5.1){$2$}; \node at (-2,-6.9){$4$};
\node at (-1.3,-5.1){$7$}; \node at (-1.3,-6.9){$6$};
\node at (-0.7,-5.1){$10$}; \node at (-0.7,-6.9){$3$};
\node at (0,-5.1){$14^*$}; \node at (0,-6.9){$17^*$};

\node at (3.3,-5.1){$22^*$}; \node at (3.3,-6.9){$19^*$};
\node at (2.7,-5.1){$21^*$}; \node at (2.7,-6.9){$16^*$};
\node at (2,-5.1){$20^*$}; \node at (2,-6.9){$15^*$};
\node at (1.3,-5.1){$12$}; \node at (1.3,-6.9){$11$};
\node at (0.7,-5.1){$18^*$}; \node at (0.7,-6.9){$9$};
\node at (4.5,-6){\vector(1,0){20}};
\node at (4.5,-6){\vector(-1,0){20}};
\foreach \i in {8,...,18} {
\node[circle, draw, fill=black, inner sep=1pt](\i) at (\i/1.5,-5.5) {};
\node[circle, draw, fill=black, inner sep=1pt](-\i) at (\i/1.5,-6.5) {};
\draw (\i) -- (-\i);
}
\node at (5.3,-5.1){$1$}; \node at (5.3,-6.9){$8$};
\node at (6,-5.1){$13^*$}; \node at (6,-6.9){$5$};
\node at (6.7,-5.1){$2$}; \node at (6.7,-6.9){$4$};
\node at (7.3,-5.1){$7$}; \node at (7.3,-6.9){$6$};
\node at (8,-5.1){$10$}; \node at (8,-6.9){$3$};
\node at (8.7,-5.1){$17^*$}; \node at (8.7,-6.9){$14^*$};
\node at (9.4,-5.1){$18^*$}; \node at (9.4,-6.9){$9$};
\node at (10,-5.1){$12$}; \node at (10,-6.9){$11$};
\node at (10.7,-5.1){$15^*$}; \node at (10.7,-6.9){$20^*$};
\node at (11.4,-5.1){$16^*$}; \node at (11.4,-6.9){$21^*$};
\node at (12.1,-5.1){$19^*$}; \node at (12.1,-6.9){$22^*$};
\node at (8.9,-3.5){\vector(0,1){20}};
\node at (9.1,-3.5){\vector(0,-1){20}};
\node at (4.2,0.2){\vector(1,0){20}};
\node at (4.2,-0.2){\vector(-1,0){20}};
\node at (4.2,0.4){$\Psi$};
\node at (8,2) {$(3,1)$};
\node[circle, draw, fill=black, inner sep=1pt, label=12] (r) at (8,1) {};
\foreach \i in {-2} {
\node[circle, draw, fill=black, inner sep=2pt](\i) at (8+\i,0) {};
\draw (r) -- (\i);
}
\foreach \i in {-1,...,0} {
\node[circle, draw, fill=black, inner sep=1pt](\i) at (8+\i,0) {};
\draw (r) -- (\i);
}
\node[circle, draw, fill=black, inner sep=1pt] at (10,0) {};
\draw (r) -- (10,0);

\node at (6,0.5) {$11$};

\node[circle, draw, inner sep=2pt] at (7,-1) {};
\draw (7,0) -- (7,-1);
\node at (7.3,0) {$2$};
\node at (7,-1.3) {$4$};

\node[circle, draw, inner sep=2pt] at (8,-1) {};
\draw (8,0) -- (8,-1);
\node at (8.3,0) {$7$};
\node at (8,-1.3) {$6$};

\foreach \i in {-1,...,2} {
\node[circle, draw, fill=black, inner sep=1pt](\i) at (10+\i,-1) {};
\draw (10,0) -- (\i);
}
\node at (10,0.3) {$1$};
\node at (9,-1.3) {$8$};
\node at (10,-1.3) {$5$};
\node at (11.3,-1) {$10$};
\node at (12,-1.3) {$9$};

\node[circle, draw, inner sep=2pt] at (11,-2) {};
\draw (11,-1) -- (11,-2);

\end{tikzpicture}
    \caption{An example of involution $\Psi$ on labelled tip-augmented plane trees}
    \label{psi}
\end{figure}
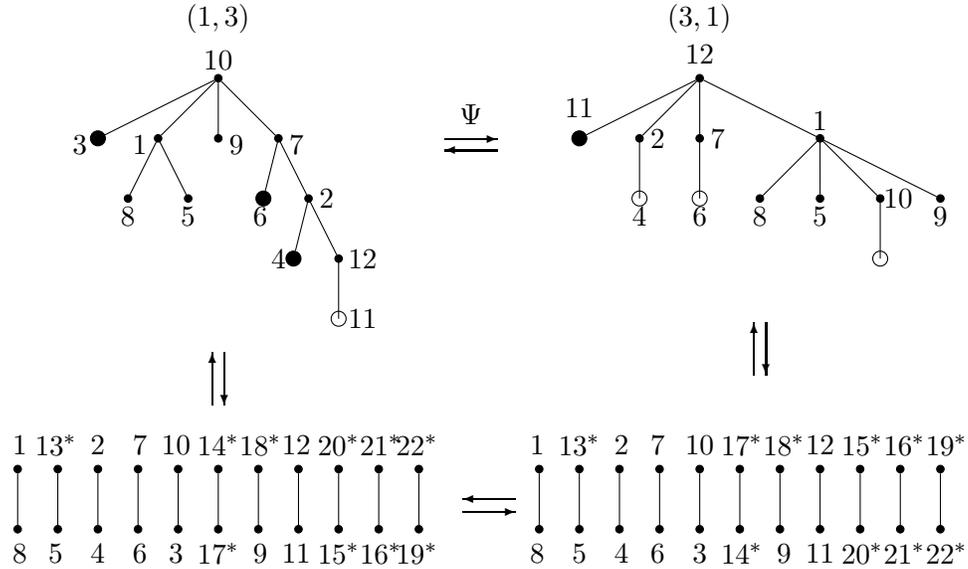

Figure~\ref{psi} illustrates how the involution $\Psi$ operates on labelled tip-augmented plane trees.

\vspace{5cm}

\noindent{\bf Acknowledgments.} The first author was partially supported by the Professional Development Program at the University of Central Florida.

\end{document}